\documentclass[11pt]{amsart}

\usepackage{amsfonts}
\usepackage{amssymb}
\usepackage{amscd}
\input xy
\xyoption{all}

\def\ZZ{{\mathbb Z}}

\def\QQ{{\mathbb Q}}
\def\PP{{\textbf P}}
\def\OO{{\mathcal O}}
\def\R{\mathbf{R}}
\def\SS{\mathcal{S}}
\def\D{\mathbf{D}}

\def\F{\mathcal{F}}
\def\E{\mathcal{E}}
\def\G{\mathcal{G}}
\def\H{\mathcal{H}}

\def\I{\mathcal{I}}
\def\J{\mathcal{J}}

\def\Pic0{{\rm Pic}^0(X)}
\def\PicY{{\rm Pic}^0(Y)}

\theoremstyle{plain}

\newtheorem{theorem}{Theorem}[section]
\newtheorem{proposition}[theorem]{Proposition}
\newtheorem{corollary}[theorem]{Corollary}
\newtheorem{lemma}[theorem]{Lemma}

\theoremstyle{definition}
\newtheorem{definition}[theorem]{Definition}
\newtheorem{remark}[theorem]{Remark}
\newtheorem{example}[theorem]{Example}

\newtheorem{conjecture/question}[theorem]{Conjecture/Question}
\newtheorem{question}[theorem]{Question}
\newtheorem{remark/definition}[theorem]{Remark/Definition}
\newtheorem{notation/terminology}[theorem]{Notation/Terminology}
\newtheorem{proposition/definition}[theorem]{Proposition/Definition}

 \theoremstyle{remark}

\begin{document}

\title{Regularity  on abelian varieties III: relationship with Generic Vanishing and applications}

\author[G. Pareschi]{Giuseppe Pareschi}
\address{Dipartamento di Matematica, Universit\`a di Roma, Tor Vergata, V.le della
Ricerca Scientifica, I-00133 Roma, Italy} \email{{\tt
pareschi@mat.uniroma2.it}}

\author[M. Popa]{Mihnea Popa}
\address{Department of Mathematics, University of Illinois at Chicago,
851 S. Morgan St., Chicago, IL 60607, USA } \email{{\tt
mpopa@math.uic.edu}}

\thanks{MP was partially supported by the NSF grant DMS 0500985,
by an AMS Centennial Fellowship, and by a Sloan Fellowship.}

\maketitle

\tableofcontents

\setlength{\parskip}{.1 in}

\markboth{G. PARESCHI and M. POPA} {\bf REGULARITY AND GENERIC
VANISHING ON ABELIAN VARIETIES}

\section{Introduction}

In previous work we have introduced the notion of
\emph{$M$-regularity} for coherent sheaves on abelian varieties
(\cite{pp1}, \cite{pp2}). This is useful because
$M$-regular sheaves enjoy strong generation
properties, in such a way that $M$-regularity on abelian varieties
presents close analogies with the classical notion of
Castelnuovo-Mumford regularity on projective spaces. Later
we studied objects in the derived
category of a smooth projective variety subject to Generic
Vanishing conditions ($GV$-\emph{objects} for short, \cite{pp4}). The main
ingredients are  Fourier-Mukai  transforms and the systematic use of
homological and commutative algebra techniques. It turns out that,
from the general perspective, $M$-regularity is a natural strenghtening of a
Generic Vanishing condition. In this paper we describe in detail the
relationship between the two notions in the case of abelian
varieties, and deduce new basic properties of both $M$-regular and
$GV$-sheaves. We also collect a few extra applications of the
generation properties of $M$-regular sheaves, mostly announced but not
contained in \cite{pp1} and \cite{pp2}. This second part of the paper is based
on our earlier preprint \cite{pp6}.

We start in \S2 by recalling  some basic definitions and results from \cite{pp4}
on $GV$-conditions, restricted to the context of the present paper (coherent sheaves
on abelian varieties). The rest of the section is
devoted to the relationship between $GV$-sheaves and $M$-regular
sheaves. More precisely, we prove a criterion, Proposition
\ref{reg-abvar}, characterizing the latter among the former:
$M$-regular sheaves are those $GV$-sheaves $\F$ for which the
Fourier-Mukai transform of the Grothendieck-dual object $\R\Delta
\F$ is a torsion-free sheaf. (This will be extended to higher regularity conditions, or 
strong Generic Vanishing conditions, in our upcoming work \cite{pp5}.)

We apply this relationship in \S3 to the basic problem of the behavior of cohomological support loci
under tensor products. We first prove that tensor products of $GV$-sheaves are again $GV$ when one of the factors is locally free, and then use
this and the torsion-freeness characterization to deduce a similar
result for $M$-regular sheaves. The question of the behavior of $M$-regularity under tensor products had been
posed to us by A. Beauville as well. It is worth mentioning that Theorem \ref{tensor} does not seem to follow
by any more standard methods.

In the other direction, in \S4 we prove a result on $GV$-sheaves based on results on $M$-regularity.
 Specifically, we show that
$GV$-sheaves on abelian varieties are nef. We deduce this from a theorem of Debarre \cite{debarre2},
stating that
$M$-regular sheaves are ample, and the results in \S2. This is especially interesting for the well-known problem of semipositivity: higher
direct images of dualizing sheaves via maps to abelian varieties are known to be $GV$ (cf. \cite{hacon}, \cite{pp3}).

In \S5 we survey generation properties of $M$-regular sheaves.
This section is mostly expository, but the presentation of some known
results, as Theorem \ref{surjectivity}$(a)\Rightarrow (b)$ (which was
proved in \cite{pp1}), is new and more natural with respect to the
Generic Vanishing perspective, providing also the new
implication $(b)\Rightarrow (a)$. In combination with well-known results of
Green-Lazarsfeld and Ein-Lazarsfeld, we deduce some basic generation properties
of the canonical bundle on a variety of maximal Albanese dimension, used in
the following section.

The second part of the paper contains miscellaneous
applications of the generation properties enjoyed by $M$-regular
sheaves on abelian varieties, extracted or reworked from our older
preprint \cite{pp6}. In \S6 we give effective results for
pluricanonical maps on irregular varieties of general type and
maximal Albanese dimension via $M$-regularity for direct images of
canonical bundles, extending work in \cite{pp1} \S5. In
particular we show, with a rather quick argument, that on a smooth
projective variety $Y$ of general type, maximal Albanese dimension,
and  whose Albanese image is not ruled by subtori, the
pluricanonical series $|3K_Y|$ is very ample outside the exceptional
locus of the Albanese map (Theorem \ref{3can}). This is a slight
strengthening, but also under a slightly stronger hypothesis, of a result
of Chen and Hacon (\cite{chenhacon}, Theorem 4.4), both statements being
generalizations of the fact that the tricanonical bundle is very ample for
curves of genus at least $2$.

In \S7.1 we look at bounding the Seshadri constant measuring the
local positivity of an ample line bundle. There is already extensive
literature on this in the case of abelian varieties (cf.
\cite{lazarsfeld1}, \cite{nakamaye}, \cite{bauer1}, \cite{bauer2},
\cite{debarre1} and also \cite{positivity} for further references).
Here we explain how the Seshadri constant of a polarization $L$ on
an abelian variety is bounded below by an asymptotic version -- and
in particular by the usual -- $M$-regularity index of the line
bundle $L$, defined in \cite{pp2} (cf. Theorem \ref{seshadri}).
Combining this with various bounds for Seshadri constants proved in
\cite{lazarsfeld1}, we obtain bounds for $M$-regularity indices
which are not apparent otherwise.

In \S7.2 we shift our attention towards a cohomological study of
Picard bundles, vector bundles on Jacobians of curves closely
related to Brill-Noether theory (cf. \cite{positivity} 6.3.C and
7.2.C for a general introduction). We combine Fourier-Mukai
techniques with the use of the Eagon-Northcott resolution for
special determinantal varieties in order to compute their regularity,
as well as that of their relatively small tensor
powers (cf. Theorem \ref{mreg-picard}). This vanishing theorem has
practical applications. In particular we recover in a more direct
fashion the main results of \cite{pp1} \S4 on the equations of the
$W_d$'s  in Jacobians, and on vanishing for pull-backs of pluritheta
line bundles to symmetric products.

By work of Mukai and others (\cite{mukai3}, \cite{mukai4},
\cite{mukai1}, \cite{umemura} and \cite{orlov}) it has emerged that
on abelian varieties the class of vector bundles most closely
resembling semistable vector bundles on curves and line bundles on
abelian varieties is that of \emph{semihomogeneous} vector bundles.
In \S7.3 we show that there exist numerical criteria for their
geometric properties like global or normal generation, based on
their Theta regularity. More generally, we give a result on the
surjectivity of the multiplication map on global sections for two
such vector bundles (cf. Theorem \ref{mixed-multiplication}). Basic
examples are the projective normality of ample line bundles on any
abelian variety, and the normal generation of the Verlinde bundles
on the Jacobian of a curve, coming from moduli spaces of vector
bundles on that curve.

\noindent \textbf{Acknowledgements.} We would like to thank Rob
Lazarsfeld for having introduced us to some of these topics and for
interesting suggestions. We also thank Christopher Hacon for
discussions, and Olivier Debarre for pointing out a mistake in
\S7.2. Finally, the second author thanks the organizers of the Clay
Workshop, Emma Previato and Montserrat Teixidor i Bigas, for
providing a few very nice days of mathematical interaction.

\section{$GV$-sheaves and $M$-regular sheaves on abelian varieties}

\noindent
{\bf $GV$-sheaves.}
We recall definitions and results from
\cite{pp4}  on  \emph{Generic Vanishing} conditions ($GV$ for
short). In relationship to the treatment of \cite{pp4} we
confine ourselves to a more limited setting, with respect to
the following three aspects: (a) we consider only
coherent sheaves (rather than complexes) subject to generic
vanishing conditions; (b) we consider only the simplest such
condition, i.e. $GV_0$, henceforth denoted $GV$; (c) we
work only on abelian varieties, with the classical Fourier-Mukai
functor associated to the Poincar\'e line bundle on $X\times \Pic0$
(rather than arbitrary integral transforms).

Let $X$ be an abelian variety of dimension $g$ over an algebraically
closed field, $\widehat X=\Pic0$, $P$ a normalized Poincar\'e bundle
on $X\times \widehat X$, and $\R \widehat \SS: \D(X) \rightarrow
\D(\widehat X)$ the standard Fourier-Mukai functor given by $\R
\widehat \SS (\F) = \R{p_{\widehat A}}_* (p_A^* \F \otimes P)$. We
denote $\R\SS: \D(\widehat X) \rightarrow \D(X)$ the functor in the
other direction defined analogously.
For a coherent sheaf $\F$ on $X$, we will consider for each $i \ge 0$ its \emph{$i$-th
cohomological support locus}
$$V^i(\F) :=\{\alpha\in \widehat X \ | \ h^i(X,\F\otimes\alpha)>0\}.$$
By base-change, the support of $R^i \widehat \SS \F$ is contained in
$V^i(\F)$.

\begin{proposition/definition}[$GV$-sheaf, \cite{pp4}]\label{def-GV}
Given a coherent sheaf $\F$ on $X$, the following  conditions are
equivalent:

\noindent
(a) ${\rm codim}~{\rm Supp}(R^i \widehat \SS \F) \ge i {\rm ~for~all~} i> 0.$

\noindent
(b) ${\rm codim}~V^i(\F) \ge i {\rm ~for~all~} i> 0.$

\noindent
If one of the above conditions is satisfied, $\F$ is called a \emph{$GV$-sheaf}.
(The proof of the equivalence is a standard base-change argument -- cf.
\cite{pp4} Lemma 3.6.)
\end{proposition/definition}

\begin{notation/terminology}
(a) ($IT_0$-sheaf). The simplest examples of $GV$-sheaves are those  such that
$V^i(\F)=\emptyset$ for every $i>0$. In this case $\F$ is said to
\emph{satisfy the Index Theorem with index 0} ($IT_0$ for
short). If $\F$ is $IT_0$ then $\R\widehat\SS \F = R^0\widehat\SS \F$, which is
a locally free sheaf.

\noindent
(b) (Weak Index Theorem).  Let $\G$ be an
object in $\D(X)$ and $k\in{\bf Z}$. $\G$ is said \emph{to satisfy
the Weak Index Theorem with index $k$} ($WIT_k$ for short), if
$R^i\widehat \SS \G =0$ for $i\ne k$. In this case we denote
$\widehat \G = R^k\widehat \SS \G $. Hence
$\R\widehat\SS \G = \widehat\G[-k]$.

\noindent
(c) The same terminology and notation holds for sheaves on $\widehat X$, or more generally objects
in $\D(\widehat X)$, considering the functor $\R\SS$.
\end{notation/terminology}

We now state a basic result from \cite{pp4} only in the special case of abelian varieties considered in this paper. In this case, with the exception of the implications from (1) to the other parts, it was
in fact proved earlier by Hacon \cite{hacon}.  We denote $\R \Delta \F : = \R \H om (\F, \OO_X)$.

\begin{theorem}\label{gv-abvar}
Let $X$ be an abelian variety and $\F$ a coherent sheaf on $X$. Then
the following are equivalent:

\noindent (1) $\F$ is a $GV$-sheaf.\hfill\break (2) For any
sufficiently positive ample line bundle $A$ on $\widehat{X}$,
$$H^i (\F \otimes \widehat{A^{-1}}) = 0, {\rm ~ for~ all~} i>0.$$
(3) $\R\Delta\F$ satisfies $WIT_g$.
\end{theorem}

\begin{proof}
This is Corollary 3.10 of \cite{pp4}, with the slight difference that
conditions (1), (2) and (3) are all stated with respect to the Poincar\'e line bundle
$P$, while condition (3)  of Corollary 3.10 of \emph{loc. cit.} holds with
respect to $P^\vee$. This can be done since, on abelian varieties,
the Poincar\' e bundle satisfies the symmetry relation $P^\vee \cong
((-1_X)\times 1_{\widehat X})^* P$. Therefore Grothendieck duality
(cf. Lemma \ref{gd} below) gives that the Fourier-Mukai functor
defined by $P^\vee$ on $X \times \widehat X$ is the same as
$(-1_X)^* \circ \R \widehat \SS$. We can also assume without loss of
generality that the ample line bundle $A$ on $\widehat X$ considered
below is symmetric.
\end{proof}

\begin{remark}
The above Theorem holds in much greater generality (\cite{pp4}, Corollary
3.10). Moreover,  in \cite{pp5} we will show that the equivalence between (1)
and (3) holds in a local setting as well. Condition (2) is a Kodaira-Kawamata-Viehweg-type
vanishing criterion. This is because, up to an \'etale cover of $X$, the vector
bundle $\widehat{A^{-1}}$ is a direct sum of copies of an ample line
bundle (cf. \cite{hacon},  and also \cite{pp4} and the
proof of Theorem \ref{nef} in the sequel).
\end{remark}

\begin{lemma}[\cite{mukai1} 3.8]\label{gd}
The Fourier-Mukai and duality functors satisfy the exchange formula:
$$\R\Delta \circ \R\widehat \SS \cong (-1_{\widehat X})^* \circ \R \widehat \SS \circ \R \Delta [g].$$
\end{lemma}

A useful immediate consequence of the equivalence of (a) and (c) of
Theorem \ref{gv-abvar}, together with Lemma \ref{gd}, is the following (cf. \cite{pp4}, Remark 3.11.):

\begin{corollary}\label{gv-duality}
If $\F$ is a $GV$-sheaf on $X$ then
$$R^i\widehat \SS \F \cong \mathcal{E}xt^i(\widehat{\R\Delta \F},\OO_{\widehat X}).$$
\end{corollary}

\noindent
{\bf $M$-regular sheaves and their characterization.}
We now recall the $M$-regularity condition, which is
simply a stronger (by one) generic vanishing condition, and relate it to the notion
of $GV$-sheaf. The reason
for the different terminology is that the notion of $M$-regularity
was discovered -- in connection with many geometric applications --
before fully appreciating its relationship with generic vanishing
theorems (see \cite{pp1}, \cite{pp2}, \cite{pp3}).

\begin{proposition/definition}
Let $\F$ be a coherent sheaf on an abelian variety $X$. The
following conditions are equivalent:

\noindent
(a) ${\rm codim}~{\rm Supp}(R^i \widehat \SS \F) > i {\rm ~for~all~} i>0.$

\noindent
(b) ${\rm codim}~V^i(\F)
> i {\rm ~for~all~} i> 0.$

\noindent
If one of the above conditions is satisfied, $\F$ is called an \emph{$M$-regular sheaf}.
\end{proposition/definition}

The proof is identical to that of Proposition/Definition \ref{def-GV}. By definition, every
$M$-regular sheaf is a $GV$-sheaf. Non-regular $GV$-sheaves are
those whose support loci have dimension as big as possible. As shown
by the next result, as a consequence of the Auslander-Buchsbaum
theorem, this is equivalent to the presence of torsion in the
Fourier transform of the Grothendieck dual object.

\begin{proposition}\label{reg-abvar}
Let $X$ be an abelian variety of dimension $g$, and let $\F$ be a
$GV$-sheaf on $X$. The following conditions are equivalent:

\noindent (1) $\F$ is M-regular.\hfill\break
(2)
$\widehat{\R\Delta\F} = \R \widehat{\mathcal{S}}(\R\Delta\F)[g]$ is
a torsion-free sheaf.\footnote{Note that it is a sheaf by Theorem
\ref{gv-abvar}.}
\end{proposition}
\begin{proof}
By Corollary  \ref{gv-duality},  $\F$ is $M$-regular if
and only if for each $i>0$
$${\rm codim} ~{\rm Supp}(\mathcal{E}xt^i(\widehat{\R\Delta \F},\OO_{\widehat X}))  > i.$$
The theorem is then a consequence of the following commutative
algebra fact, which is surely known to the experts.
\end{proof}

\begin{lemma}\label{torsion}
Let $\G$ be a coherent sheaf on a
smooth variety $X$. Then $\G$ is torsion-free if and only if
${\rm codim} ~{\rm Supp}(\mathcal{E}xt^i(\G,\OO_X)) > i$ for all $i>0$.
\end{lemma}

\begin{proof}
If $\G$ is torsion free then it is a subsheaf of a  locally free
sheaf $\E$. From the exact sequence
$$0\longrightarrow \G \rightarrow \E
\longrightarrow \E/\G\longrightarrow 0$$
it follows that, for $i > 0$,
$\mathcal{E}xt^i(\G,\OO_{\widehat X}) \cong
\mathcal{E}xt^{i+1}(\E/\G,\OO_{\widehat X})$. But then a well-known
consequence of the Auslander-Buchsbaum Theorem applied to $\E/\G$
implies that
$${\rm codim} ~{\rm Supp}(\mathcal{E}xt^i(\G,\OO_{\widehat X})) > i,
{\rm ~for~all~}i >0.$$

Conversely, since $X$ is smooth, the functor
$\R\mathcal{H}om(\,\cdot\,,\OO_X)$ is an involution on $\D(X)$.
Thus there is a spectral sequence
$$E^{ij}_2 := \mathcal{E}xt^i \Bigl((\mathcal{E}xt^j (\G, \OO_X),\OO_X\Bigr)
\Rightarrow H^{i -j} = \mathcal{H}^{i-j} \G  = \begin{cases}\G
&\hbox{if $i= j$}\cr 0&\hbox{otherwise}\cr\end{cases}.$$ If $ {\rm
codim} ~{\rm Supp}(\mathcal{E}xt^i(\G,\OO_X)) > i $ $
{\rm ~for~all~}i > 0,$ then  $\mathcal{E}xt^i \Bigl(\mathcal{E}xt^j
(\G, \OO_X),\OO_X \Bigr) = 0$ for all $i,j$ such that  $j > 0$ and
$i - j \le 0$, so the only $E^{ii}_{\infty}$ term which might be
non-zero is $E^{00}_{\infty}$. But the differentials coming into
$E^{00}_p$ are always zero, so we get a sequence of inclusions
$$\F=H^0= E^{00}_{\infty} \subset\ldots \subset E^{00}_3 \subset E^{00}_2.$$
The extremes give precisely the injectivity of the natural map
 $\G \rightarrow \G^{**}$. Hence $\G$ is torsion free.
\end{proof}

\begin{remark}
It is worth noting that in the previous proof, the
fact that we are working on an abelian varieties is of no importance.
In fact, an extension of Proposition \ref{reg-abvar} holds in the generality of
\cite{pp4}, and even in a local setting,  as it will be shown in \cite{pp5}.
\end{remark}

\section{Tensor products of $GV$ and $M$-regular sheaves}

We now address the issue of preservation of bounds on the
codimension of support loci under tensor products. Our main result
in this direction  is (2) of Theorem \ref{tensor} below, namely that
the tensor product of two $M$-regular sheaves on an abelian variety
is $M$-regular, provided that one of them is locally free. Note that
the same result holds for Castelnuovo-Mumford regularity on
projective spaces (\cite{positivity}, Proposition 1.8.9). We do not
know whether the same holds if one removes the local freeness
condition on $\E$ (in the case of Castelnuovo-Mumford regularity it
does not).

Unlike the previous section, the proof of the result is quite
specific to abelian varieties. One of the essential
ingredients is  Mukai's main inversion result (cf. \cite{mukai1},
Theorem 2.2), which states that the functor $\R \widehat \SS$ is an
equivalence of derived categories and, more precisely,
\begin{equation}\label{inversion}
\mathbf{R}\mathcal{S}\circ\mathbf{R}\widehat{\mathcal{S}}\cong
(-1_A)^{*} [-g] ~{\rm ~and~}~ \mathbf{R}\widehat{\mathcal{S}}
\circ\mathbf{R}\mathcal{S}\cong (-1_{\widehat{A}})^{*} [-g].
\end{equation}
Besides this, the argument uses  the
characterization of $M$-regularity among $GV$-sheaves given by
Proposition \ref{reg-abvar}.

\begin{proposition}\label{ab-vanishing2}
Let $\F$ be a $GV$-sheaf and $H$ a locally free sheaf satisfying $IT_0$
on an abelian variety $X$. Then $\F\otimes H$ satisfies
$IT_0$.
\end{proposition}
\begin{proof}
Consider any $\alpha\in {\rm Pic}^0(X)$. Note that $H\otimes \alpha$
also satisfies $IT_0$, so $\R \widehat{\SS} (H \otimes \alpha) = R^0
\widehat{\SS} (H \otimes \alpha)$ is a vector bundle $N_{\alpha}$ on
$\widehat X$. By Mukai's inversion theorem (\ref{inversion})
 $N_\alpha$ satisfies $WIT_g$ with respect to $\R \SS$ and $H\otimes \alpha \cong \R\SS ( (-1_X)^*N_{\alpha})
[g]$. Consequently for all $i$ we have
\begin{equation}\label{bah}
H^i (X, \F \otimes H\otimes
\alpha) \cong H^i (X, \F \otimes \R\SS ( (-1_{\widehat
X})^*N_{\alpha}) [g]).
\end{equation}
But a basic exchange formula
for integral transforms (\cite{pp4}, Lemma 2.1) states, in the
present context, that
\begin{equation}\label{exchange}
H^i (X, \F \otimes \R\SS ( (-1_{\widehat X})^*N_{\alpha})
[g])\cong H^i (Y, \R\widehat\SS\F \underline\otimes (-1_{\widehat X})^*N_{\alpha}
[g]).
\end{equation}
Putting (\ref{bah}) and (\ref{exchange}) together, we get that
\begin{equation}\label{final}
H^i (X, \F \otimes H\otimes \alpha) \cong  H^i
(Y, \R\widehat\SS\F \underline\otimes (-1_{\widehat X})^*N_{\alpha}
[g])= H^{g+i}(Y, \R\widehat\SS\F \underline\otimes (-1_{\widehat X})^*N_{\alpha}).
\end{equation}
The hypercohomology groups on the right hand side are computed by the
spectral sequence
$$E_2^{jk}: = H^j(Y, R^k\SS\F \otimes (-1_{\widehat
X})^*N_{\alpha} [g])\Rightarrow H^{j+k}(Y, \R\widehat\SS\F
\underline\otimes (-1_{\widehat X})^*N_{\alpha} [g]).$$ Since $\F$
is $GV$, we have the vanishing of $H^j(Y, R^k\SS\F\otimes
(-1_{\widehat X})^*N_{\alpha} [g])$ for $j+k> g$, and from this it
follows that the hypercohomology groups in (\ref{final}) are zero
for $i>0$.
\end{proof}

\begin{theorem}\label{tensor}
Let $X$ be an abelian variety, and $\F$ and $\E$ two coherent
sheaves on $X$, with $\E$ locally free.
\newline
\noindent (1) If $\F$ and $\E$ are $GV$-sheaves,
 then $\F\otimes \E$ is a $GV$-sheaf.
\newline
\noindent (2) If $\F$ and $\E$ are $M$-regular,
then $\F\otimes \E$ is $M$-regular.
\end{theorem}
\begin{proof}
(1) Let $A$ be a sufficiently ample line bundle on $\widehat{X}$.
Then, by Theorem \ref{gv-abvar}(2), $\E\otimes \widehat{A^{-1}}$
satisfies $IT_0$. By Proposition \ref{ab-vanishing2}, this implies
that $(\F \otimes \E) \otimes \widehat{A^{-1}}$ also satisfies
$IT_0$. Applying Theorem \ref{gv-abvar}(2) again, we deduce that
$\F\otimes \E$ is $GV$.

\noindent (2) Both $\F$ and $\E$ are $GV$, so (1) implies that $\F\otimes \E$ is also a $GV$.
We use Proposition \ref{reg-abvar}. This implies to begin with that $\R
\widehat \SS(\R\Delta\F)$ and $\R \widehat \SS(\R \Delta \E) \cong
\R \widehat \SS( \E^\vee)$ are torsion-free sheaves (we harmlessly
forget about what degree they live in). Going backwards, it also
implies that we are done if we show that $\R \widehat \SS(\R
\Delta(\F\otimes \E))$ is torsion free. But note that
$$\R \widehat \SS(\R \Delta(\F\otimes \E)) \cong
\R \widehat \SS(\R \Delta \F\otimes \E^\vee)\cong \R \widehat \SS(\R
\Delta \F) \underline * \R \widehat \SS(\E^\vee)$$ where $\underline
*$ denotes the (derived) Pontrjagin product of sheaves on abelian varieties,
and the last isomorphism is the exchange of Pontrjagin and tensor
products under the Fourier-Mukai functor (cf. \cite{mukai1} (3.7)).
Note that this derived Pontrjagin product is in fact an honest
Pontrjagin product, as we know that all the objects above are
sheaves. Recall that by definition the Pontrjagin product of two
sheaves $\G$ and $\H$ is simply $\G * \H := m_* (p_1^* \G \otimes
p_2^* \H)$, where \hfill\break $m : \widehat X\times \widehat X
\rightarrow \widehat X$ is the group law on $\widehat X$. Since $m$
is a surjective morphism, if $\G$ and $\H$ are torsion-free, then so
is $p_1^* \G \otimes p_2^* \H$ and its push-forward $\G * \H$.
\end{proof}

\begin{remark}
As mentioned  in \S2, Generic Vanishing
conditions can be naturally defined for objects in the derived
category, rather than sheaves (see \cite{pp4}). In this more general
setting,  (1) of Theorem \ref{tensor} holds for
$\F\underline\otimes\G$, where $\F$ is any $GV$-\emph{object} and $\E$
any $GV$-\emph{sheaf}, while (2) holds for $\F$ any $M$-regular
\emph{object} and $\E$ any $M$-regular \emph{locally free sheaf}.
The proof is the same.
\end{remark}

\section{Nefness of $GV$-sheaves}

Debarre has shown in \cite{debarre2} that every $M$-regular sheaf on an abelian variety
is ample. We deduce from this and Theorem \ref{gv-abvar} that $GV$-sheaves satisfy the analogous weak positivity.

\begin{theorem}\label{nef}
Every $GV$-sheaf on an abelian variety is nef.
\end{theorem}
\begin{proof}
\emph{Step 1.} We first reduce to the case when the abelian variety
$X$ is principally polarized. For this, consider $A$ any ample line
bundle on $\widehat X$. By Theorem \ref{gv-abvar} we know that the
$GV$-condition is equivalent to the vanishing
$$H^i (\F \otimes \widehat{A^{-m}}) = 0, ~{\rm for ~all~} i> 0,  ~{\rm and~ all~} m >> 0.$$
But $A$ is the pullback ${\hat\psi}^*L$ of a principal polarization
$L$ via an isogeny $\hat\psi: \widehat X \rightarrow \widehat Y$
(cf. \cite{lange} Proposition 4.1.2). We then have
$$ 0 = H^i (\F \otimes \widehat{A^{-m}}) \cong H^i (\F \otimes \widehat{({\hat\psi}^*(L^{-m}))})
\cong H^i (\F \otimes {\psi}_* \widehat{M^{-m}})\cong H^i ({\psi}^*
\F \otimes \widehat{M^{-m}}).$$ Here ${\psi}$ denotes the dual
isogeny. (The only thing that needs an explanation is the next to
last isomorphism, which is the commutation of the Fourier-Mukai
functor with isogenies, \cite{mukai1} 3.4.) But this implies that
${\psi}^* \F$ is also $GV$, and since nefness is preserved by
isogenies this completes the reduction step.

\noindent \emph{Step 2.} Assume now that $X$ is principally
polarized by $\Theta$. As above, we know that
$$H^i (\F \otimes \widehat{\OO(-m\Theta)}\otimes \alpha) = 0, ~{\rm for ~all~} i> 0, {\rm~all~} \alpha \in {\rm Pic}^0(X)
 ~{\rm and~ all~} m >> 0.$$
If we denote by $\phi_m: X \rightarrow X$ multiplication by $m$,
i.e. the isogeny induced by $m\Theta$, then this implies that
$$H^i (\phi_m ^* \F \otimes \OO(m\Theta) \otimes \beta) = 0, ~{\rm for ~all~} i> 0 {\rm~and~all~}
\beta \in {\rm Pic}^0(X)$$ as $\phi_m^*  \widehat{\OO(-m\Theta)}
\cong \bigoplus \OO(m\Theta)$ by \cite{mukai1} Proposition 3.11(1).
This means that the sheaf $\phi_m ^* \F \otimes \OO(m\Theta)$
satisfies $IT_0$ on $X$, so in particular it is $M$-regular. By
Debarre's result \cite{debarre2} Corollary 3.2, it is then ample.

But $\phi_m$ is a finite cover, and $\phi_m^* \Theta \equiv
m^2\Theta$. The statement above is then same as saying that, in the
terminology of \cite{positivity} \S6.2, the
$\QQ$-twisted\footnote{Note that the twist is indeed only up to
numerical equivalence.} sheaf $\F<\frac{1}{m}\cdot \Theta>$ on $X$
is ample, since $\phi_m^* (\F<\frac{1}{m}\cdot \Theta>)$ is an
honest ample sheaf. As $m$ goes to $\infty$, we see that $\F$ is a
limit of ample $\QQ$-twisted sheaves, and so it is nef by
\cite{positivity} Proposition 6.2.11.
\end{proof}

Combining the result above with the fact that higher direct images
of canonical bundles are $GV$ (cf. \cite{pp3} Theorem 5.9), we
obtain the following result, one well-known instance of which is
that the canonical bundle of any smooth subvariety of an abelian
variety is nef.

\begin{corollary}\label{nef-image}
Let $X$ be a smooth projective variety and $a: Y \rightarrow X$ a (not necessarily
surjective) morphism to an abelian variety. Then $R^j a_* \omega_Y$ is a nef
sheaf on $X$ for all $j$.
\end{corollary}

One example of an immediate application of Corollary  is to
integrate a result of Peternell-Sommese in the general picture.

\begin{corollary}[\cite{ps}, Theorem 1.17]\label{covering}
Let $a: Y \rightarrow X$ be a finite surjective morphism of smooth
projective varieties, with $X$ an abelian variety. Then the vector
bundle $E_a$
 is nef.
\end{corollary}
\begin{proof}
By duality we have $a_* \omega_Y \cong \OO_X \oplus E_a$. Thus $E_a$
is a quotient of $a_* \omega_Y$, so by Corollary \ref{nef-image} it
is nef.
\end{proof}

\section{Generation properties of $M$-regular sheaves on abelian varieties}

The interest in the notion of $M$-regularity comes from the fact
that $M$-regular sheaves on abelian varieties have strong generation
properties. In this respect, $M$-regularity on abelian varieties
parallels the notion of Castelnuovo-Mumford regularity on projective
spaces (cf. the survey \cite{pp3}). In this section we survey the
basic results about generation properties of $M$-regular sheaves.
The presentation is somewhat new, since the proof of the basic
result (the implication $(a)\Rightarrow (b)$ of Theorem
\ref{surjectivity} below) makes use of the relationship between
$M$-regularity and $GV$-sheaves (Proposition \ref{reg-abvar}). The
argument in this setting turns out to be more natural, and provides
as a byproduct the reverse implication $(b)\Rightarrow (a)$, which
is new. 

\noindent {\bf Another characterization of $M$-regularity.}
$M$-regular sheaves on abelian varieties are characterized as
follows:

\begin{theorem}(\cite{pp1}, Theorem 2.5)\label{surjectivity}
Let $\F$ be a $GV$-sheaf on an abelian variety $X$. Then the
following conditions are equivalent:

\noindent
(a) $\F$ is $M$-regular.

\noindent
(b) For every locally free sheaf $H$ on $X$
satisfying $IT_0$, and for every  non-empty Zariski open set
$U\subset \widehat X$, the sum of multiplication maps of global sections
$${\mathcal M}_U:\bigoplus_{\alpha\in U}H^0(X,\F\otimes \alpha)\otimes
H^0(X,H\otimes \alpha^{-1})
\buildrel{\oplus m_\alpha}\over \longrightarrow H^0(X,\F\otimes H)$$ is
surjective.
\end{theorem}

\begin{proof}
Since $\F$ is a $GV$-sheaf, by Theorem \ref{gv-abvar} the transform of
$\R\Delta\F$ is a sheaf in degree $g$, i.e.
$\R\widehat\SS(\R\Delta\F)=\widehat{\R\Delta\F}[-g]$. If $H$ is a
coherent sheaf satisfying $IT_0$ then $\R\widehat\SS
H=\widehat{H}$, a locally free sheaf in degree $0$.
It turns out that the following natural map is an isomorphism
\begin{equation}\label{iso} {\rm Ext}^g(H,\R\Delta \F)\buildrel{\sim}\over\longrightarrow {\rm
Hom}(\widehat{H},\widehat{\R\Delta\F}).
\end{equation}
This simply follows from Mukai's Theorem $(\ref{inversion})$, which yields that $${\rm
Ext}^g(H,\R\Delta F)={\rm Hom}_{\D(X)}(H,\R\Delta\F[g])\cong{\rm
Hom}_{\D(\widehat X)} (\widehat{H},\widehat{\R\Delta\F})={\rm
Hom}(\widehat{H},\widehat{\R\Delta\F}).$$

\noindent \emph{Proof of (a) $\Rightarrow$ (b).} Since
$\widehat{\R\Delta\F}$ is torsion-free by Proposition
\ref{reg-abvar}, the evaluation map at the fibres
\begin{equation}\label{eval}  {\rm
Hom}(\widehat{H},\widehat{\R\Delta\F})\rightarrow\prod_{\alpha\in U}
{\mathcal
Hom}(\widehat{H},\widehat{\R\Delta\F})\otimes_{\OO_{\widehat
X,\alpha}} k(\alpha)
\end{equation}
is injective for all open sets $U\subset\Pic0$. Therefore, composing
with the isomorphism (\ref{iso}), we get an injection
\begin{equation}\label{compo}{\rm Ext}^g(H,\R\Delta
F)\rightarrow\prod_{\alpha\in U}
{\mathcal
Hom}(\widehat{H},\widehat{\R\Delta\F})\otimes_{\OO_{\widehat X,\alpha}}
k(\alpha).
\end{equation}
By base-change, this is the dual map of the map in (b), which is
therefore surjective.

\noindent \emph{Proof of (b) $\Rightarrow$ (a).} Let $A$ be an ample
symmetric line bundle on $\widehat X$. From Mukai's Theorem
$(\ref{inversion})$, it follows that $A^{-1} =\widehat{H_A}$, where
$H_A$ is a locally free sheaf on $X$ satisfying $IT_0$ and such that
$\widehat{H_A}=A^{-1}$. We have that (b) is equivalent to the
injectivity of (\ref{compo}). We now take $H=H_A$ in both 
(\ref{iso}) and (\ref{compo}). The facts that (\ref{iso}) is an
isomorphism and that (\ref{compo}) is injective yield
 the
injectivity, \emph{for all open sets $U\subset\Pic0$}, of the
evaluation map at fibers
$$H^0(\widehat{\R\Delta\F}\otimes A)\overset{ev_U}{\longrightarrow}\prod_{\alpha\in U}
(\widehat{\R\Delta\F}\otimes A)\otimes_{\OO_{\widehat X ,\alpha}}
k(\alpha).$$ Letting $A$ be sufficiently positive so that
$\widehat{\R\Delta\F}\otimes A$ is globally generated, this is
equivalent to the torsion-freeness of
$\widehat{\R\Delta\F}$\footnote{Note that the kernel of $ev_U$
generates a torsion subsheaf of $\widehat{\R\Delta\F}\otimes A$
whose support is contained in the complement of $U$.} and hence, by
Proposition \ref{reg-abvar}, to the $M$-regularity of $\F$.
\end{proof}

\noindent
{\bf Continuous global generation and global generation.}
Recall first the following:

\begin{definition}[\cite{pp1}, Definition 2.10]\label{continuous}
Let $Y$ be a variety equipped with a morphism $a:Y\rightarrow X$ to
an abelian variety $X$.

\noindent
(a) A sheaf $\F$ on $Y$ is \emph{continuously globally generated with respect to $a$}
if the sum of evaluation maps
$${\mathcal Ev}_U:\bigoplus_{\alpha\in U} H^0(\F\otimes a^*\alpha)\otimes a^*\alpha^{-1} \longrightarrow \F$$
is surjective  for every non-empty open subset $U\subset \Pic0$.

\noindent (b) More generally, let  $T$ be a proper subvariety of
$Y$. The sheaf $\F$ is said to be \emph{continuously globally
generated with respect to $a$ away from $T$} if  \ ${\rm Supp}({\rm
Coker} ~{\mathcal E}v_U)\subset T$ \  for every non-empty open
subset $U\subset \Pic0$.

\noindent
(c) When $a$ is the Albanese
morphism, we will suppress $a$ from the terminology, speaking of
\emph{continuously globally generated} (resp. \emph{continuously globally generated
away from $T$}) sheaves.
\end{definition}\label{CGG}

In Theorem \ref{surjectivity}, taking $H$ to be a sufficiently
positive line bundle on $X$ easily yields  (cf \cite{pp1},
Proposition 2.13):

\begin{corollary}\label{M-reg}
An $M$-regular sheaf on $X$ is continuously globally generated.
\end{corollary}

\noindent
The relationship between continuous global generation and global
generation comes from:

\begin{proposition}[\cite{pp1}, Proposition 2.12]\label{GG}
(a) In the setting of Definition \ref{CGG}, let $\F$ (resp. $A$)
be a coherent sheaf on $Y$ (resp. a line
bundle, possibly supported on a subvariety $Z$ of $Y$), both
continuously globally generated. Then $\F\otimes A\otimes a^*\alpha$
is globally generated for all $\alpha\in\Pic0$.

\noindent
(b) More generally, let $\F$ and $A$ as above. Assume
that $\F$ is continuously globally generated away from $T$ and that
$A$ is continuously globally generated away from $W$.  Then $F\otimes
A\otimes a^*\alpha$ is globally generated away from $T\cup W$ for all
$\alpha\in\Pic0$.
\end{proposition}

The proposition is proved via the classical method of
\emph{reducible sections}, i.e. those sections of the form
$s_\alpha\cdot t_{-\alpha}$, where $s_\alpha$ (resp. $t_{-\alpha}$)
belongs to $H^0(F\otimes a^*\alpha)$ (resp. $H^0(A\otimes
a^*\alpha^{-1})$).

\noindent
{\bf Generation properties on varieties
of maximal Albanese dimension via Generic Vanishing.}
The above results give effective generation criteria
once one has effective \emph{Generic Vanishing} criteria ensuring
that the dimension of the cohomological support loci is not too
big. The main example of such a criterion is the Green-Lazarsfeld
Generic Vanishing Theorem for the canonical line bundle of an
irregular variety,  proved in \cite{greenlaz1} and
further refined in \cite{greenlaz2} using the deformation theory of
cohomology groups.\footnote{More recently, Hacon \cite{hacon} has given a different proof,
based on the Fourier-Mukai transform and Kodaira
Vanishing.  Building in part on Hacon's ideas, several
extensions of this result are given in \cite{pp4}.} For the purposes of this
paper, it is enough to state the Generic Vanishing
Theorem in the case of varieties $Y$ of \emph{maximal Albanese
dimension}, i.e. such that the Albanese map $a: Y \rightarrow {\rm Alb}(Y)$
is generically finite onto its image. More generally, we consider a morphism
$a:Y\rightarrow X$ to an abelian variety $X$. Then,  as in
\S1 one can consider the \emph{cohomological support loci}
$V_a^i(\omega_Y)=\{\alpha\in \Pic0 \>|\>h^i(\omega_Y\otimes
a^*\alpha)>0\>\}$. (In case $a$ is the Albanese map we will suppress
$a$ from the notation.)

The result of Green-Lazarsfeld (see also \cite{einlaz} Remark 1.6)
states that, \emph{if the morphism $a$ is generically finite, then}
$${\rm codim~}V_a^i(\omega_Y) \ge i ~ {\rm ~for~all~} i> 0.$$
Moreover, in \cite{greenlaz2} it is proved
that \emph{$V_a^i(\omega_Y)$ are unions of translates of subtori}.
Finally, an argument of Ein-Lazarsfeld \cite{einlaz} yields that,
\emph{if there exists an $i>0$ such that ${\rm
codim}~V_a^i(\omega_Y)=i$, then the image of $a$ is ruled by subtori
of $X$.} All of this implies the following typical application of the
concept of $M$-regularity.

\begin{proposition}\label{key-canonical}
Assume that $\dim Y=\dim a(Y)$ and that $a(Y)$ is not ruled by tori.
Let $Z$  be the exceptional locus of $a$, i.e. the inverse image via
$a$ of the locus of points in $a(Y)$ having non-finite fiber. Then:
\newline\noindent (i) $a_*\omega_Y$ is an M-regular sheaf on $X$.
\newline\noindent (ii) $a_*\omega_Y$ is continuously globally
generated.
\newline\noindent (iii) $\omega_Y$ is continuously globally generated away from
$Z$.
\newline\noindent (iv) For all  $k\ge 2$, $\omega_Y^{\otimes k}\otimes a^*\alpha$ is globally generated away from $Z$ for any $\alpha\in \Pic0$.
\end{proposition}
\begin{proof}
By Grauert-Riemenschneider vanishing, $R^ia_*\omega_Y=0$ for all
$i\ne 0$. By the Projection Formula we get
$V^i_a(\omega_Y)=V^i(a_*\omega_Y)$. Combined with the Ein-Lazarsfeld
result, (i) follows. Part (ii) follows from Corollary \ref{M-reg}.
For (iii) note that, as with global generation (and by a similar
argument), continuous global generation is preserved by finite maps:
if $a$ is finite and $a_* \F$ is continuously globally generated,
then $\F$ is continuously globally generated. (iv) for $k=2$ follows
from (iii) and Proposition \ref{GG}. For arbitrary $k\ge 2$ it
follows in the same way by induction (note that if a sheaf $\F$ is
such that $\F\otimes a^*\alpha$ is globally generated away from $Z$
for every $\alpha\in \Pic0$, then it is continuously globally
generated away from $Z$).
\end{proof}

\section{Pluricanonical maps of irregular varieties of maximal
Albanese dimension}

One of the most elementary results about projective
embeddings is that every curve of general type can be embedded in
projective space by the tricanonical line bundle. This is sharp
for curves of genus two. It turns out that this result can be
generalized to arbitrary dimension, namely to varieties
of maximal Albanese dimension. In fact, using Vanishing and
Generic Vanishing Theorems and the Fourier-Mukai transform, Chen and
Hacon proved that \emph{for every smooth
complex variety of general type and maximal Albanese dimension $Y$
such that $\chi(\omega_Y)>0$, the tricanonical line bundle
$\omega_Y^{\otimes 3}$ gives a birational map}  (cf. \cite{chenhacon}, Theorem 4.4).
The main point of this section is that
the concept of $M$-regularity (combined of course with vanishing results) provides a quick and
conceptually simple proof of on one hand a slightly more explicit version of
the Chen-Hacon Theorem, but  on the other hand under a slightly more restrictive hypothesis. We show the following:

\begin{theorem}\label{3can}
Let $Y$ be a smooth projective complex  variety of general type and
maximal Albanese dimension. If the Albanese image of $Y$ is not
ruled by tori, then $\omega_Y^{\otimes 3}$ is very ample away from the
exceptional locus of the Albanese map.
\end{theorem}

Here the \emph{exceptional locus} of the Albanese map $a:Y\rightarrow
{\rm Alb}(Y)$ is $Z=a^{-1}(T)$, where $T$ is the locus of points
in ${\rm Alb}(Y)$ over which the fiber of $a$ has positive dimension.

\begin{remark}\label{hypothesis}
A word about the hypothesis of the Chen-Hacon Theorem and of Theorem
\ref{3can} is in order. As a consequence of the Green-Lazarsfeld
Generic Vanishing Theorem (end of \S3), it follows that
$\chi(\omega_Y)\ge 0$ for every variety $Y$ of maximal Albanese
dimension. Moreover, Ein-Lazarsfeld \cite{einlaz} prove that for $Y$
of maximal Albanese dimension, \emph{if $\chi(\omega_Y)=0$, then
$a(Y)$ is ruled by subtori of ${\rm Alb}(Y)$}. In dimension $\ge 3$
there exist examples of varieties of general type and maximal
Albanese dimension with $\chi(\omega_Y)=0$ (cf. \emph{loc. cit.}).
\end{remark}

In the course of the proof we will invoke $\J(Y, \parallel
L\parallel)$, the \emph{asymptotic multiplier ideal sheaf} associated
to a complete linear series $|L|$ (cf. \cite{positivity} \S11). One knows
that, given a line bundle $L$ of non-negative Iitaka dimension,
\begin{equation}\label{a-i}
H^0(Y, L \otimes \J(\parallel L\parallel)) = H^0(Y, L),
\end{equation}
i.e. the zero locus of $\J(\parallel L\parallel)$ is contained in
the base locus of $|L|$ (\cite{positivity}, Proposition 11.2.10). Another
basic property we will use is that, for every $k$,
\begin{equation}\label{sequence}
\J(\parallel L^{\otimes (k+1)}\parallel)\subseteq\J(\parallel L^{\otimes k} \parallel).
\end{equation}
(Cf. \cite{positivity}, Theorem 11.1.8.)  A first standard result is

\begin{lemma}\label{invariance}
Let $Y$ be a smooth projective complex variety of general
type. Then:

\noindent
(a) $h^0(\omega^{\otimes m}_Y\otimes \alpha)$ is constant for all $\alpha\in
\PicY$ and for all $m>1$.

\noindent (b) The zero locus of $\J(\parallel \omega_Y^{\otimes
(m-1)}\parallel)$ is contained in the base locus of
$\omega_Y^{\otimes m}\otimes\alpha$, for all $\alpha\in \PicY$.
\end{lemma}
\begin{proof}
Since bigness is a numerical property, all line bundles $\omega_Y\otimes\alpha$
are big, for $\alpha\in\PicY$. By  Nadel Vanishing for
asymptotic multiplier ideals (\cite{positivity}, Theorem 11.2.12)
$$H^i(Y, \omega^{\otimes m}_Y\otimes \beta \otimes \J(\parallel(\omega_Y\otimes\alpha)^{\otimes (m-1)}\parallel))=0$$
for all $i>0$ and all $\alpha,\beta\in \Pic0$. Therefore, by the invariance of the
Euler characteristic,
$$h^0(Y,\omega^{\otimes m}_Y\otimes \beta\otimes \J(\parallel(\omega_Y\otimes\alpha)^{\otimes (m-1)}\parallel))={\rm constant}=\lambda_\alpha$$
for all $\beta\in \PicY$.
Now
$$h^0(Y,\omega^{\otimes m}_Y\otimes \beta\otimes \J(\parallel(\omega_Y\otimes\alpha)^{\otimes (m-1)}\parallel)) \le h^0(Y,\omega_Y^{m}\otimes\beta)$$
for all $\beta\in\Pic0$ and, because of (\ref{a-i}) and
(\ref{sequence}), equality holds for $\beta=\alpha^{m}$. By
semicontinuity it follows that $h^0(Y,\omega_Y^{\otimes
m}\otimes\beta)=\lambda_\alpha$ for all $\beta$
contained in a Zariski open set $U_\alpha$ of $\Pic0$ which contains
$\alpha^m$. Since this is true for all $\alpha$, the statement
follows. Part (b) follows from the previous argument.
\end{proof}

\begin{lemma}\label{asymptotic}
Let $Y$ be a smooth projective complex
variety of general type and maximal Albanese dimension, such that its Albanese image
is not ruled by tori. Let $Z$
be the exceptional locus of its Albanese map. Then, for every
$\alpha\in\PicY$:

\noindent (a) the zero-locus of
$\J(\parallel\omega_Y\otimes\alpha\parallel)$ is contained
(set-theoretically) in $Z$.

\noindent
(b) $\omega_Y^{\otimes 2}\otimes
\alpha\otimes \J(\parallel\omega_Y\parallel)$ is globally generated away from $Z$.
\end{lemma}
\begin{proof}
(a) By (\ref{a-i}) and (\ref{sequence}) the zero locus of
$\J(\parallel\omega\otimes\alpha\parallel)$ is contained in the base
locus of $\omega^{\otimes 2}\otimes \alpha^2$. By Proposition
\ref{key-canonical}, the base locus of $\omega^{\otimes
2}\otimes\alpha^2$ is contained $Z$. (b) Again by Proposition
\ref{key-canonical}, the base locus of $\omega^{\otimes
2}\otimes\alpha$ is contained in $Z$.  By Lemma \ref{invariance}(b),
the zero locus of $\J(\parallel\omega_Y\parallel)$ is contained in
$Z$.
\end{proof}

\begin{proof}\emph{(of Theorem \ref{3can})}
As above, let $a: Y \rightarrow {\rm Alb}(Y)$ be the Albanese map and let $Z$ be the exceptional locus of $a$.
 As in the proof of Prop. \ref{key-canonical},  the Ein-Lazarsfeld result at the end of \S3 (see also Remark \ref{hypothesis}), the hypothesis implies that $a_* \omega_Y$ is $M$-regular, so $\omega_Y$ is continuously globally generated away from $Z$. We make the following:

\noindent
{\bf Claim.} \emph{For every
$y\in Y-Z$, the sheaf $a_*(I_y \otimes
\omega_Y^{\otimes 2} \otimes\J(\parallel\omega_Y\parallel))$ is $M$-regular.}

\noindent We first see how the Claim implies Theorem \ref{3can}. The
statement of the Theorem is equivalent to the fact that, \emph{ for
any $y\in Y-Z$, the sheaf $I_y\otimes \omega_Y^{\otimes 3}$ is
globally generated away from $Z$. } By Corollary  \ref{M-reg}, the
Claim yields that $a_*(I_y \otimes \omega_Y^{\otimes 2}
\otimes\J(\parallel\omega_Y\parallel))$ is continuously globally
generated. Therefore $I_y \otimes  \omega_Y^{\otimes 2} \otimes
\J(\parallel\omega_Y\parallel)$ is continuously globally generated
away from $Z$. Hence, by Proposition \ref{GG}, $I_y \otimes
\omega_Y^{\otimes 3}\otimes \J(\parallel\omega_Y\parallel)$ is
globally generated away from $Z$. Since the zero locus of
$\J(\parallel\omega_Y\parallel)$ is contained in $Z$ (by Lemma
\ref{asymptotic})(a)), the Theorem follows from the Claim.

\noindent
\emph{Proof of the Claim.}
We consider the standard exact sequence
\begin{equation}\label{standard}
0\rightarrow I_y\otimes \omega^{\otimes 2}_Y\otimes
\alpha\otimes  \J(\parallel\omega_Y\parallel)
\rightarrow \omega^{\otimes 2}_Y\otimes \alpha  \otimes \J(\parallel\omega_Y\parallel)
\rightarrow {(\omega^{\otimes 2}_Y\otimes \alpha  \otimes \J(\parallel\omega_Y\parallel))}_{|y} \rightarrow 0.
\end{equation}
(Note that $y$ does not lie in the zero locus of
$\J(\parallel\omega_Y\parallel)$.) By Nadel Vanishing for asymptotic
multiplier ideals, $H^i(Y, \omega^{\otimes 2}_Y\otimes \alpha\otimes
\J(\parallel\omega_Y\parallel))=0$ for all $i>0$ and $\alpha\in
\PicY$. Since, by Lemma \ref{asymptotic}, $y$ is not in the base
locus of $\omega^{\otimes 2}_Y\otimes \alpha\otimes
\J(\parallel\omega_Y\parallel)$, taking cohomology in
(\ref{standard}) it follows that
\begin{equation}\label{vanishing}H^i(Y, I_y\otimes\omega^{\otimes
2}_Y\otimes \alpha\otimes \J(\parallel\omega_Y\parallel))=0
\end{equation}
for all $i>0$ and $\alpha\in \Pic0$ as well.
Since $y$ does not belong to the exceptional locus of $a$, the map
$a_*( \omega^{2}_Y\otimes \J(\parallel\omega_Y\parallel))\rightarrow
a_*({(\omega^{2}_Y\otimes \J(\parallel\omega_Y\parallel))}_{|y})$ is
surjective.  On the other hand, since $a$ is generically finite, by
a well-known extension of Grauert-Riemenschneider vanishing,\break
$R^ia_*( \omega^{\otimes 2}_Y\otimes \J(\parallel\omega_Y\parallel))$
vanishes for all $i>0$.\footnote{The proof of this is identical to that of the usual 
Grauert-Riemenschneider vanishing theorem in \cite{positivity} \S4.3.B, replacing 
Kawamata-Viehweg vanishing with Nadel vanishing.}  
Therefore (\ref{standard}) implies also that for all $i>0$
\begin{equation}\label{relative-vanishing}
R^ia_*(I_y\otimes \omega^{\otimes 2}_Y\otimes \J(\parallel\omega_Y\parallel))=0.
\end{equation}
Combining (\ref{vanishing}) and (\ref{relative-vanishing}) one gets,
by  projection Formula, that the sheaf $a_*(I_y\otimes
\omega^{\otimes 2}_Y\otimes \J(\parallel\omega_Y\parallel))$ is
$IT_0$ on $X$, hence $M$-regular.
\end{proof}

\begin{remark} It follows from the proof that
$\omega_Y^{\otimes 3}\otimes\alpha$ is very ample away from $Z$ for all
$\alpha\in\PicY$ as well.
\end{remark}

\begin{remark}[The Chen-Hacon Theorem]
The reader might wonder why, according to the above quoted theorem
of Chen-Hacon, the tricanonical bundle of varieties of general type
and maximal Albanese dimension is birational (but not necessarily
very ample outside the Albanese exceptional locus) even under the
weaker assumption that $\chi(\omega_Y)$ is positive, which does not
ensure the continuous global generation of $a_*\omega_Y$. The point
is that, according to Generic Vanishing, if the Albanese dimension
is maximal, then $\chi(\omega_Y)>0$ implies
$h^0(\omega_Y\otimes\alpha)>0$ for all $\alpha\in\PicY$. Hence, even
if $\omega_Y$ is not necessarily continuously globally generated
away of some subvariety of $Y$, the following condition holds:  for
\emph{general} $y\in Y$, there is a Zariski open set $U_y\subset
\PicY$ such that $y$ is not a base point of $\omega_Y\otimes\alpha$
for all $\alpha\in U_y$. Using the same argument of Proposition
\ref{GG} -- based on reducible sections -- it follows that such $y$
is not a base point of $\omega_Y^{\otimes 2}\otimes\alpha$ for
\emph{all} $\alpha\in\PicY$. Then the Chen-Hacon Theorem follows by
an argument analogous to that of Theorem \ref{3can}.
\end{remark}

To complete the picture, it remains to analyze the case of varieties
$Y$ of maximal Albanese dimension and $\chi(\omega_Y)=0$. Chen and Hacon
prove that if the Albanese dimension is maximal, then $\omega_Y^{\otimes 6}$ is
always birational (and $\omega_Y^{\otimes 6}\otimes\alpha$ as well). The same
result can be made slightly more precise as follows, extending also results in \cite{pp1} \S5:

\begin{theorem}
If $Y$ is a smooth projective complex variety of maximal
Albanese dimension then, for all $\alpha\in\PicY$,
$\omega_Y^{\otimes 6}\otimes\alpha$ is very ample away from the exceptional
locus of the Albanese map. Moreover, if $L$ a big line bundle on $Y$, then
$(\omega_Y \otimes L)^{\otimes 3}\otimes\alpha$ gives a birational map.
\end{theorem}

The proof is similar to that of Theorem \ref{3can}, and left to the
interested reader. For example, for the first part the point is that, by Nadel Vanishing
for asymptotic multiplier ideals, $H^i(Y, \omega_Y^{\otimes 2}\otimes\alpha\otimes
\J(\parallel\omega_Y\parallel))=0$ for all
$\alpha\in\PicY$. Hence, by the same argument using
Grauert-Riemenschneider vanishing, $a_*(\omega_Y^{\otimes 2}\otimes
\J(\parallel\omega_Y\parallel))$ is $M$-regular.

Finally, we remark that in \cite{chenhacon}, Chen-Hacon also prove
effective birationality results for pluricanonical maps of
irregular varieties of arbitrary Albanese dimension (in function of
the minimal power for which the corresponding pluricanonical map on
the general Albanese fiber is birational). It is likely that the methods above apply
to this context as well.

\section{Further applications of $M$-regularity}

\subsection{$M$-regularity indices and Seshadri constants.}

Here we express a natural relationship between Seshadri constants of
ample line bundles on abelian varieties and the $M$-regularity
indices of those line bundles as defined in \cite{pp2}. This result
is a theoretical improvement of the lower bound for Seshadri
constants proved in \cite{nakamaye}. In the opposite direction,
combined with the results of \cite{lazarsfeld1}, it provides bounds
for controlling $M$-regularity. For a general overview of Seshadri
constants, in particular the statments used below, one can consult
\cite{positivity} Ch.I \S5.

We start by recalling the basic definition from \cite{pp2} and by
also looking at a slight variation. We will denote by $X$ an abelian
variety of dimension $g$ over an algebraically closed field and by
$L$ an ample line bundle on $X$.

\begin{definition}
The $M$-\emph{regularity index} of $L$ is defined as
$$m(L):={\rm max}\{l~|~L\otimes
m_{x_1}^{k_1}\otimes \ldots \otimes m_{x_p}^{k_p}~{\rm is ~}M{\rm
-regular~ for ~all~distinct~}$$
$$x_1,\ldots,x_p\in X {\rm~with~} \Sigma k_i=l\}.$$
\end{definition}

\begin{definition}
We also define a related invariant, associated to just one given
point $x\in X$:
$$p(L,x) := {\rm max}\{l~|~L\otimes m_x^l~{\rm is ~}M{\rm -regular}\}.$$
The definition does not depend on $x$ because of the homogeneity of
$X$, so we will denote this invariant simply by $p(L)$.
\end{definition}

Our main interest will be in the asymptotic versions of these
indices, which turn out to be related to the Seshadri constant
associated to $L$.

\begin{definition}
The \emph{asymptotic $M$-regularity index} of $L$ and its punctual
counterpart are defined as
$$\rho(L) := \underset{n}{{\rm sup}} \frac{m(L^n)}{n} {\rm~~and~~}
\rho^\prime (L) := \underset{n}{{\rm sup}} \frac{p(L^n)}{n}.$$
\end{definition}

\noindent The main result of this section is:

\begin{theorem}\label{seshadri}
We have the following inequalities:
$$\epsilon(L) = \rho^{\prime}(L) \geq \rho(L)\geq 1.$$
In particular $\epsilon(L)\ge {\rm max}\{m(L),1\}$.
\end{theorem}

This improves a result of Nakamaye (cf. \cite{nakamaye} and the
references therein). Nakamaye also shows that $\epsilon(L) = 1$ for
some line bundle $L$ if and only if $X$ is the product of an
elliptic curve with another abelian variety. As explained in
\cite{pp2} \S3, the value of $m(L)$ is reflected in the geometry of
the map to projective space given by $L$. Here is a basic example:

\begin{example}
If $L$ is very ample -- or more generally gives a birational
morphism outside a codimension $2$ subset --  then $m(L)\ge 2$, so
by the Theorem above $\epsilon (L)\ge 2$. Note that on an arbitrary
smooth projective variety very ampleness implies in general only
that $\epsilon(L,x)\ge 1$ at each point.
\end{example}

The proof of Theorem \ref{seshadri} is a simple application Corollary \ref{M-reg}
and Proposition \ref{GG}, via the results of
\cite{pp2} \S3. We use the relationship with the notions of $k$-jet
ampleness and separation of jets. Denote by $s(L,x)$ the largest
number $s\ge 0$ such that $L$ separates $s$-jets at $x$. Recall also
the following:

\begin{definition}
A line bundle $L$ is called $k$-\emph{jet ample}, $k\geq 0$, if the
restriction map
$$H^0(L)\longrightarrow H^0(L\otimes \OO_X/
m_{x_1}^{k_1}\otimes \ldots \otimes m_{x_p}^{k_p})$$ is surjective
for any distinct points $x_1,\ldots,x_p$ on $X$ such that $\Sigma
k_i =k+1$. Note that if $L$ is $k$-jet ample, then it separates
$k$-jets at every point.
\end{definition}

\begin{proposition}[\cite{pp2} Theorem 3.8 and Proposition 3.5]\label{both_ways}
(i) $L^n$ is ($n + m(L) -2$)-jet ample, so in particular
$s(L^n,x)\ge n + m(L) -2$.
\newline
\noindent (ii) If $L$ is $k$-jet ample, then $m(L)\geq k+1$.
\end{proposition}

This points in the direction of local positivity, since one way to
interpret the Seshadri constant of $L$ is (independently of $x$):
$$\epsilon (L) = \underset{n}{{\rm sup}} \frac{s(L^n,x)}{n}.$$
To establish the connection with the asymptotic invariants above we
also need the following:

\begin{lemma}\label{comparison}
For any $n\geq 1$ and any $x\in X$ we have $s(L^{n+1},x)\geq
m(L^n)$.
\end{lemma}
\begin{proof}
This follows immediately from Corollary \ref{M-reg} and Proposition \ref{GG}:
if $L^n\otimes m_{x_1}^{k_1}\otimes \ldots \otimes m_{x_p}^{k_p}$ is $M$-regular,
then $L^{n+1}\otimes m_{x_1}^{k_1}\otimes \ldots \otimes
m_{x_p}^{k_p}$ is globally generated, and so by \cite{pp2} Lemma
3.3, $L^{n+1}$ is $m(L)$-jet ample.
\end{proof}

\begin{proof}(\emph{of Theorem \ref{seshadri}.})
Note first that for every $p\ge 1$ we have
\begin{equation}\label{scaling}
m(L^n)\ge m(L) + n -1,
\end{equation}
which follows immediately from the two parts of Proposition
\ref{both_ways}. In particular $m(L^n)$ is always at least $n-1$,
and so $\rho(L)\ge 1$. Putting together the definitions,
(\ref{scaling}) and Lemma \ref{comparison}, we obtain the main
inequality $\epsilon(L) \geq \rho(L)$. Finally, the asymptotic
punctual index computes precisely the Seshadri constant. Indeed, by
completely similar arguments as above, we have that for any ample
line bundle $L$ and any $p\geq 1$ one has
$$p(L^n)\ge s(L^n,x)~{\rm and~} s(L^{n+1},x)\ge p(L^n, x).$$
The statement follows then from the definition.
\end{proof}

\begin{remark}
What the proof above shows is that one can give an interpretation
for $\rho(L)$ similar to that for $\epsilon (L)$ in terms of
separation of jets. In fact $\rho(L)$ is precisely the ``asymptotic
jet ampleness" of $L$, namely:
$$\rho (L) = \underset{n}{{\rm sup}} \frac{a(L^n)}{n},$$
where $a(M)$ is the largest integer $k$ for which a line bundle $M$
is $k$-jet ample.
\end{remark}

\begin{question}
Do we always have $\epsilon(L) = \rho(L)$? Can one give independent
lower bounds for $\rho(L)$ or $\rho^\prime(L)$ (which would then
bound Seshadri constants from below)?
\end{question}

In the other direction, there are numerous bounds on Seshadri
constants, which in turn give bounds for the $M$-regularity indices
that (at least to us) are not obvious from the definition. All of
the results in \cite{positivity} Ch.I \S5 gives some sort of bound.
Let's just give a couple of examples:

\begin{corollary}
If $(J(C), \Theta)$ is a principally polarized Jacobian, then
$m(n\Theta)\le \sqrt{g}\cdot n$. On an arbitrary abelian variety,
for any principal polarization $\Theta$ we have $m(n\Theta)\le
(g!)^{\frac{1}{g}}\cdot n$.
\end{corollary}
\begin{proof}
It is shown in \cite{lazarsfeld1} that $\epsilon(\Theta)\le
\sqrt{g}$. We then apply Theorem \ref{seshadri}. For the other bound
we use the usual elementary upper bound for Seshadri constants,
namely $\epsilon(\Theta) \leq (g!)^{\frac{1}{g}}$.
\end{proof}

\begin{corollary}
If $(A, \Theta)$ is a very general PPAV, then there exists at least
one $n$ such that $p(n\Theta)\ge \frac{2^{\frac{1}{g}}}{4}
(g!)^{\frac{1}{g}}\cdot n$.
\end{corollary}
\begin{proof}
Here we use the lower bound given in \cite{lazarsfeld1} via a result
of Buser-Sarnak.
\end{proof}

There exist more specific results on $\epsilon (\Theta)$ for
Jacobians (cf. \cite{debarre1}, Theorem 7), each giving a
corresponding result for $m(n\Theta)$. We can ask however:

\begin{question}
Can we calculate $m(n\Theta)$ individually on Jacobians, at least
for small $n$, in terms of the geometry of the curve?
\end{question}

\begin{example}[Elliptic curves]
As a simple example, the question above has a clear answer for
elliptic curves. We know that on an elliptic curve $E$ a line bundle
$L$ is $M$-regular if and only if ${\rm deg}(L)\ge 1$, i.e. if and
only if $L$ is ample. From the definition of $M$-regularity we see
then that if ${\rm deg}(L)=d>0$, then $m(L)=d-1$. This implies that
on an elliptic curve $m(n\Theta) = n-1$ for all $n\ge 1$. This is
misleading in higher genus however; in the simplest case we have the
following general statement: \emph{If $(X,\Theta)$ is an irreducible
principally polarized abelian variety of dimension at least $2$,
then $m(2\Theta)\geq 2$}. This is an immediate consequence of the
properties of the Kummer map. The linear series $|2\Theta|$ induces a
$2:1$ map of $X$ onto its image in $\PP^{2^g-1}$, with injective
differential. Thus the cohomological support locus for
$\OO(2\Theta)\otimes m_x\otimes m_y$ consists of a finite number of points, while
the one for $\OO(2\Theta)\otimes m_x^2$ is empty.
\end{example}

\subsection{Regularity of Picard bundles and vanishing on symmetric products.}

In this subsection we study the regularity of Picard bundles
over the Jacobian of a curve, twisted by positive multiples of the
theta divisor. Some applications to the degrees of equations
cutting out special subvarieties of Jacobians are drawn in the second
part. Let $C$ be a smooth curve of genus $g\geq 2$, and denote by
$J(C)$ the Jacobian of $C$. The objects we are interested in are the
Picard bundles on $J(C)$:  a line bundle $L$ on $C$ of degree
$n\geq 2g-1$ -- seen as a sheaf on $J(C)$ via an Abel-Jacobi
embedding of $C$ into $J(C)$ -- satisfies $IT_0$, and the
Fourier-Mukai transform $E_L=\widehat L$ is called an \emph{ n-th
Picard bundle} . When possible, we omit the dependence on $L$
and write simply $E$. Note that any other such $n$-th Picard bundle
$E_M$, with $M\in {\rm Pic}^n(C)$, is a translate of $E_L$.
The line bundle $L$ induces an identification between $J(C)$ and ${\rm
Pic}^n(C)$, so that the projectivization of $E$ -- seen as a vector
bundle over ${\rm Pic}^n(C)$ -- is the symmetric product $C_n$ (cf.
\cite{acgh} Ch.VII \S2).

The following theorem is the main cohomological result we are aiming
for. It is worth noting that Picard bundles are known to be negative
(i.e with ample dual bundle), so vanishing theorems are not
automatic. To be very precise, everything that follows holds if $n$
is assumed to be at least $4g-4$. (However the value of $n$ does not
affect the applications.)

\begin{theorem}\label{mreg-picard}
For every $1\leq k\leq g-1$, $\otimes^k E\otimes\OO(\Theta)$
satisfies $IT_0$.
\end{theorem}

\noindent
Before proving the Theorem, we record the following preliminary:

\begin{lemma}\label{fm-picard}
For any $k\geq 1$, let $\pi_k: C^k \rightarrow J(C)$ a desymmetrized
Abel-Jacobi mapping and let $L$ be a line bundle on $C$ of degree
$n>>0$ as above. Then ${\pi_k}_*(L\boxtimes \ldots\boxtimes L)$
satisfies $IT_0$, and
$$({\pi_k}_*(L\boxtimes \ldots\boxtimes L))^{\widehat{}}= \otimes^k E,$$
where $E$ is the $n$-th Picard bundle of $C$.
\end{lemma}
\begin{proof}
The first assertion is clear. Concerning the second
assertion note that, by definition, ${\pi_k}_*(L\boxtimes
\ldots\boxtimes L)$ is the Pontrjagin product $L*\ldots *L$. By the
exchange of Pontrjagin and tensor product under the Fourier-Mukai
transform (\cite{mukai1} (3.7)), it follows that $(L*\ldots
*L)^{\widehat{}} \cong \widehat L\otimes\ldots\otimes \widehat
L=\otimes^kE$.
\end{proof}

\begin{proof} \emph{(of Theorem \ref{mreg-picard})}\footnote{We
are grateful to Olivier Debarre for pointing out a numerical mistake
in the statement, in a previous version of this paper.} We will use
loosely the notation $\Theta$ for any translate of the canonical
theta divisor. The statement of the theorem becomes then equivalent
to the vanishing
$$h^i(\otimes^k E\otimes \OO(\Theta))=0,~\forall~ i>0, ~\forall~ 1\leq k\leq g-1.$$
To prove this vanishing we use the Fourier-Mukai transform. The
first point is that Lemma \ref{fm-picard} above, combined with
Grothendieck duality (Theorem \ref{gd} above), tells us precisely
that $\otimes^k E$ satisfies $WIT_g$, and, by Mukai inversion
theorem (Theorem \ref{inversion}) its Fourier transform is
$$\widehat{\otimes^k E}=(-1_J)^*{\pi_k}_* (L\boxtimes\ldots\boxtimes L).$$
Using, once again, the fact that the Fouerier-Mukai tranform is an
equivalence, we have the following sequence of isomorphisms:
$$H^i(\otimes^k E\otimes \OO(\Theta))\cong {\rm Ext}^i(\OO(-\Theta), \otimes^k E)
\cong {\rm Ext}^i(\widehat{\OO(-\Theta)}, \widehat{\otimes^k E})$$
$$\cong {\rm Ext}^i(\OO(\Theta), (-1_J)^*{\pi_k}_*(L\boxtimes\ldots\boxtimes L))\cong
H^i((-1_J)^*{\pi_k}_*(L\boxtimes\ldots\boxtimes L)\otimes
\OO(-\Theta))$$ (here we are
 the fact that both $\OO(-\Theta)$ and
$\otimes^k E$ satisfy $WIT_g$ and that
$\widehat{\OO(-\Theta)}=\OO(\Theta)$).

As we are loosely writing $\Theta$ for any translate, multiplication
by $-1$ does not influence the vanishing, so the result follows if
we show:
$$h^i({\pi_k}_*(L\boxtimes\ldots\boxtimes L)\otimes \OO(-\Theta))=0, ~\forall~i>0.$$
Now the image $W_k$  of the Abel-Jacobi map $u_k:C_k \rightarrow
J(C)$ has rational singularities (cf. \cite{kempf2}), so we only
need to prove the vanishing:
$$h^i(u_k^*({\pi_k}_*(L\boxtimes\ldots\boxtimes L)\otimes \OO(-\Theta)))=0, ~\forall~i>0.$$
Thus we are interested in the skew-symmetric part of the cohomology
group $H^i(C^k, (L\boxtimes\ldots\boxtimes L)\otimes
\pi_k^*\OO(-\Theta))$, or, by Serre duality that of
$$H^i(C^k, ((\omega_C\otimes L^{-1})\boxtimes\ldots\boxtimes (\omega_C\otimes L^{-1}))
\otimes \pi_k^*\OO(\Theta)), ~{\rm for}~i<k.$$ At this stage we can
essentially invoke a Serre vanishing type argument, but it is worth
noting that the computation can be in fact made very concrete. For
the identifications used next we refer to \cite{izadi} Appendix 3.1.
As $k\le g-1$, we can write
$$\pi_k^*\OO(\Theta)\cong ((\omega_C\otimes A^{-1})\boxtimes\ldots\boxtimes
(\omega_C\otimes A^{-1}))\otimes \OO(-\Delta),$$ where $\Delta$ is
the union of all the diagonal divisors in $C^k$ and $A$ is a line
bundle of degree $g-k-1$. Then the skew-symmetric part of the
cohomology groups we are looking at is isomorphic to
$$S^i H^1(C, \omega_C^{\otimes 2}\otimes A^{-1}\otimes L^{-1})  \otimes \wedge^{k-i}
H^0(C, \omega_C^{\otimes 2}\otimes A^{-1}\otimes L^{-1}),$$ and
since for $1\leq k\leq g-1$ and $n\ge 4g-4$ the degree of the line
bundle  $\omega_C^{\otimes 2}\otimes A^{-1}\otimes L^{-1}$ is
negative, this vanishes precisely for $i<k$.
\end{proof}

An interesting consequence of the vanishing result for Picard
bundles proved above is a new -- and in some sense more classical --
way to deduce Theorem 4.1 of \cite{pp1} on the M-regularity  of
twists of ideal sheaves $\I_{W_d}$ on the Jacobian $J(C)$. This
theorem has a number of applications to the equations of the $W_d$'s
inside $J(C)$, and also to vanishing results for pull-backs  of
theta divisors to symmetric products. For this circle of ideas we
refer the reader to \cite{pp1} \S4. For any $1\leq d\leq g-1$, $g
\geq 3$, consider $u_d :C_d\longrightarrow J(C)$ to be an
Abel-Jacobi mapping of the symmetric product (depending on the
choice of a line bundle of degree $d$ on $C$), and denote by $W_d$
the image of $u_d$ in $J(C)$.

\begin{theorem}\label{two-theta}
For every $1\leq d\leq g-1$, $\I_{W_d}(2\Theta)$ satisfies $IT_0$.
\end{theorem}
\begin{proof}
We have to prove that:
$$h^i(\I_{W_d}\otimes \OO(2\Theta)\otimes \alpha) = 0,~\forall~i>0,~\forall~\alpha\in
{\rm Pic}^0(J(C)).$$ In the rest of the proof, by $\Theta$ we will
understand generically any translate of the canonical theta divisor,
and so $\alpha$ will disappear from the notation.

It is well known that $W_d$ has a natural determinantal structure,
and its ideal is resolved by an Eagon-Northcott complex. We will
chase the vanishing along this complex. This setup is precisely the
one used by Fulton and Lazarsfeld in order to prove for example the
existence theorem in Brill-Noether theory --  for explicit details
on this cf. \cite{acgh} Ch.VII \S2. Concretely, $W_d$ is the
"highest" degeneracy locus of a map of vector bundles
$$\gamma: E\longrightarrow F,$$
where ${\rm rk}F=m$ and ${\rm rk}E=n=m+d-g+1$, with $m>>0$
arbitrary. The bundles $E$ and $F$ are well understood: $E$ is the
$n$-th Picard bundle of $C$, discussed above, and $F$ is a direct
sum of topologically trivial line bundles. (For simplicity we are
again moving the whole construction on $J(C)$ via the choice of a
line bundle of degree $n$.) In other words, $W_d$ is scheme
theoretically the locus where the dual map
$$\gamma^*:F^*\longrightarrow E^*$$
fails to be surjective. This locus is resolved by an Eagon-Northcott
complex (cf. \cite{kempf}) of the form:
$$0\rightarrow \wedge^m F^*\otimes S^{m-n}E\otimes {\rm det}E\rightarrow \ldots
\rightarrow \wedge^{n+1} F^*\otimes E\otimes {\rm det}E\rightarrow
\wedge^n F^* \rightarrow \I_{W_d}\rightarrow 0.$$ As it is known
that the determinant of $E$ is nothing but $\OO(-\Theta)$, and since
$F$ is a direct sum of topologically trivial line bundles, the
statement of the theorem follows by chopping this into short exact
sequences, as long as we prove:
$$h^i(S^k E\otimes \OO(\Theta))=0,~\forall~ i>0, ~\forall~ 1\leq k\leq m-n=g-d-1.$$
Since we are in characteristic zero, $S^k E$ is naturally a direct
summand in $\otimes^k E$, and so it is sufficient to prove that:
$$h^i(\otimes^k E\otimes \OO(\Theta))=0,~\forall~ i>0, ~\forall~ 1\leq k\leq g-d-1.$$
But this follows from Theorem \ref{mreg-picard}.
\end{proof}

\begin{remark}
Using \cite{pp1} Proposition 2.9, it follows that
$\I_{W_d}(k\Theta)$ satisfies $IT_0$ for all $k\ge 2$.
\end{remark}

\begin{remark}
It is conjectured, based on a connection with minimal cohomology classes (cf. \cite{pp7} for a discussion), that the only nondegenerate subvarieties $Y$ of a principally polarized abelian variety
$(A, \Theta)$  such that $\I_Y (2\Theta)$ satisfies $IT_0$ are precisely the $W_d$'s above, in Jacobians, and the Fano surface of lines in the intermediate Jacobian of the cubic threefold.
\end{remark}

\begin{question}
What is the minimal $k$ such that
 $\I_{W_d^r}(k\Theta)$ is $M$-regular, for $r$ and $d$ arbitrary?
\end{question}

We describe below one case in which the answer can already be given,
namely that of the singular locus of the Riemann theta divisor on a
non-hyperelliptic jacobian. It should be noted that in this case we
do not have that  ${\mathcal I}_{W^1_{g-1}}(2\Theta)$ satisfies $IT_0$
any more (but rather ${\mathcal I}_{W^1_{g-1}}(3\Theta)$ does,
by the same \cite{pp1} Proposition 2.9).

\begin{proposition}\label{singtheta}
${\mathcal I}_{W^1_{g-1}}(2\Theta)$ is $M$-regular.
\end{proposition}
\begin{proof} It follows from the results of \cite{vgi}
that $$ h^i({\mathcal I}_{W^1_{g-1}}\otimes {\mathcal
O}(2\Theta)\otimes\alpha)=\begin{cases}0&\hbox{for $i\ge g-2$,
$\forall\alpha\in {\rm Pic}^0(J(C))$}\cr 0&\hbox{for $0<i<g-2$,
$\forall\alpha\in {\rm Pic}^0(J(C))$ such that $\alpha\ne {\mathcal
O}_{J(C)}.$}\cr
\end{cases}$$
For the reader's convenience, let us briefly recall the relevant
points from Section 7 of \cite{vgi}. We denote for simplicity, via
translation, $\Theta=W_{g-1}$, (so that $W^1_{g-1}= {\rm
Sing}(\Theta$)). In the first place, from the exact sequence
$$0\rightarrow {\mathcal O}(2\Theta)\otimes\alpha\otimes \OO(-\Theta)\rightarrow {\mathcal
I}_{W^1_{g-1}}(2\Theta)\otimes\alpha\rightarrow {\mathcal
I}_{W^1_{g-1}/\Theta}(2\Theta)\otimes\alpha\rightarrow 0$$ it
follows that
$$h^i( J(C) ,{\mathcal
I}_{W^1_{g-1}}(2\Theta)\otimes\alpha)=h^i(\Theta, {\mathcal
I}_{W^1_{g-1}/\Theta}(2\Theta)\otimes\alpha)\hbox{ for $i>0$.}$$
Hence one is reduced to a computation on $\Theta$.  It is a standard
fact  (see e.g. \cite{vgi}, 7.2) that, via the Abel-Jacobi map
$u=u_{g-1}:C_{g-1}\rightarrow \Theta\subset J(C)$,
$$h^i(\Theta, {\mathcal
I}_{W^1_{g-1}/\Theta}(2\Theta)\otimes\alpha)=h^i(C_{g-1}, L^{\otimes
2}\otimes \beta\otimes {\mathcal I}_Z),$$ where
$Z=u^{-1}(W^1_{g-1})$, $L=u^*{\mathcal O}_X(\Theta)$ and
$\beta=u^*\alpha$. We now use the standard exact sequence
(\cite{acgh}, p.258):
$$0\rightarrow T_{C_{g-1}}\buildrel{du}\over\rightarrow
H^1(C,{\mathcal O}_C)\otimes{\mathcal O}_{C_{g-1}}\rightarrow
L\otimes {\mathcal I}_Z\rightarrow 0.$$ Tensoring with $L\otimes
\beta$, we see that it is sufficient to prove that
$$H^i(C_{g-1},T_{C_{g-1}}\otimes L\otimes\beta )=0, \> \forall i\ge
2, \> \forall \beta\ne {\mathcal O}_{C_{g-1}}.$$ To this end we use
the well known fact (cf. \emph{loc. cit.}) that
$$T_{C_{g-1}}\cong p_*{\mathcal O}_D(D)$$
where $D\subset C_{g-1}\times C$ is the universal divisor and $p$ is
the projection onto the first factor. As $p_{|D}$ is finite, the
degeneration of the Leray spectral sequence and the projection
formula ensure that
$$h^i(C_{g-1}, T_{C_{g-1}}\otimes L\otimes\beta)=h^i(D, {\mathcal
O}_D(D)\otimes p^*(L\otimes \beta)),$$ which are zero for $i\ge 2$
and $\beta$ non-trivial by \cite{vgi}, Lemma 7.24.
\end{proof}

\subsection{Numerical study of semihomogeneous vector bundles.}

An idea that originated in work of Mukai is that on abelian
varieties the class of vector bundles to which the theory of line
bundles should generalize naturally is that of
\emph{semihomogeneous} bundles (cf. \cite{mukai1}, \cite{mukai3},
\cite{mukai4}). These vector bundles are semistable, behave nicely
under isogenies and Fourier transforms, and have a Mumford type
theta group theory as in the case of line bundles (cf.
\cite{umemura}). The purpose of this section is to show that this
analogy can be extended to include effective global generation and
normal generation statements dictated by specific numerical
invariants measuring positivity. Recall that normal generation is
Mumford's terminology for the surjectivity of the multiplication map
$H^0(E) \otimes H^0(E)\rightarrow H^0(E^{\otimes 2})$.

In order to set up a criterion for normal generation, it is useful
to introduce the following notion, which parallels the notion of
Castelnuovo-Mumford regularity.

\begin{definition}\label{theta-reg}
A coherent sheaf $\F$ on a polarized abelian variety $(X, \Theta)$
is called $m$-$\Theta$-\emph{regular} if $\F((m-1)\Theta)$ is
$M$-regular.
\end{definition}

The relationship with normal generation comes from (3) of the following
``abelian" Castelnuovo-Mumford Lemma. Note that (1) is Corollary \ref{M-reg} plus Proposition \ref{GG}.

\begin{theorem}[\cite{pp1}, Theorem 6.3]\label{acm}
Let $\F$ be a $0$-$\Theta$-regular coherent sheaf on $X$. Then:
\newline
\noindent (1) $\F$ is globally generated.
\newline
\noindent (2) $\F$ is $m$-$\Theta$-regular for any $m\geq 1$.
\newline
\noindent (3) The multiplication map
$$H^0(\F(\Theta))\otimes H^0(\OO(k\Theta))\longrightarrow H^0(\F((k+1)\Theta))$$
is surjective for any $k\geq 2$.
\end{theorem}

\noindent {\bf Basics on semihomogeneous bundles.} Let $X$ be an
abelian variety of dimension $g$ over an algebraically closed field.
As a general convention, for a numerical class $\alpha$ we will use
the notation $\alpha >0$ to express the fact that $\alpha$ is ample.
If the class is represented by an effective divisor, then the
condition of being ample is equivalent to $\alpha^g >0$. For a line
bundle $L$ on $X$, we denote by $\phi_L$ the isogeny defined by $L$:
$$\begin{array}{cccc}
\phi_{L}:& X& \longrightarrow & {\rm Pic}^{0}(X)\cong\widehat{X}\\
&x &\leadsto& t_{x}^*L\otimes L^{-1}.
\end{array}
$$

\begin{definition}(\cite{mukai3})
A vector bundle $E$ on $X$ is called \emph{semihomogeneous} if for
every $x\in X$, $t_x^*E\cong E\otimes \alpha$, for some $\alpha\in
{\rm Pic}^0(X)$.
\end{definition}

Mukai shows in \cite{mukai3} \S6 that the semihomogeneous bundles
are Gieseker semistable (while the simple ones -- i.e. with no
nontrivial automorphisms -- are in fact stable). Moreover, any
semihomogeneous bundle has a Jordan-H\" older filtration in a strong
sense.

\begin{proposition}(\cite{mukai3} Proposition 6.18)
Let $E$ be a semihomogeneous bundle on $X$, and let $\delta$ be the
equivalence class of $\frac{{\rm det}(E)}{{\rm rk}(E)}$ in
$NS(X)\otimes_\ZZ \QQ$. Then there exist simple semihomogeneous
bundles $F_1, \ldots, F_n$ whose corresponding class is the same
$\delta$, and semihomogeneous bundles $E_1, \ldots, E_n$,
satisfying:
\begin{itemize}
\item $E\cong \bigoplus_{i=1}^n E_i$.
\smallskip
\item Each $E_i$ has a filtration whose factors are all isomorphic to $F_i$.
\end{itemize}
\end{proposition}

Since the positivity of $E$ is carried through to the factors of a
Jordan-H\"older filtration as in the Proposition above, standard
inductive arguments allow us to immediately reduce the study below
to the case of simple semihomogeneous bundles, which we do freely in
what follows.

\begin{lemma}\label{sh-properties}
Let $E$ be a simple semihomogeneous bundle of rank $r$ on $X$.

\noindent (1) (\cite{mukai3}, Proposition 7.3) There exists an
isogeny $\pi: Y\rightarrow X$ and a line bundle $M$ on $Y$ such that
$\pi^*E\cong \underset{r}{\bigoplus}M$.

\noindent (2) (\cite{mukai3}, Theorem 5.8(iv)) There exists an
isogeny $\phi: Z\rightarrow X$ and a line bundle $L$ on $Z$ such
that $\phi_* L=E$.
\end{lemma}

\begin{lemma}\label{index-theorem}
Let $E$ be a nondegenerate (i.e. $\chi(E)\neq 0$) simple
semihomogeneous bundle on $X$. Then exactly one cohomology group
$H^i(E)$ is nonzero, i.e. $E$ satisfies the Index Theorem.
\end{lemma}
\begin{proof}
This follows immediately from the similar property of the line
bundle $L$ in Lemma \ref{sh-properties}(2).
\end{proof}

\begin{lemma}\label{strong}
A semihomogeneous bundle $E$ is $m$-$\Theta$-regular if and only if
$E((m-1)\Theta)$ satisfies $IT_0$.
\end{lemma}
\begin{proof}
The more general fact that an $M$-regular semihomogeneous bundle
satisfies $IT_0$ follows quickly from Lemma \ref{sh-properties}(1)
above. More precisely the line bundle $M$ in its statement is forced
to be ample since it has a twist with global sections and positive
Euler characteristic.
\end{proof}

\noindent {\bf A numerical criterion for normal generation.} The
main result of this section is that the normal generation of a
semihomogeneous vector bundle is dictated by an explicit numerical
criterion. We assume all throughout that all the semihomogeneous
vector bundles involved satisfy the minimal positivity condition,
namely that they are $0$-$\Theta$-regular, which in particular is a
criterion for global generation by Theorem \ref{acm}.
We will in fact prove a criterion which
guarantees the surjectivity of multiplication maps for two arbitrary
semihomogeneous bundles. This could be seen as an analogue of
Butler's theorem \cite{butler} for semistable bundles on
curves.

\begin{theorem}\label{mixed-multiplication}
Let $E$ and $F$ be semihomogeneous bundles on $(X,\Theta)$, both
$0$-$\Theta$-regular. Then the multiplication maps
$$H^0(E)\otimes H^0(t_x^*F)\longrightarrow H^0(E\otimes t_x^*F)$$
are surjective for all $x\in X$ if  the following holds:
$$\frac{1}{r_F}\cdot c_1(F(-\Theta)) + \frac{1}{r_E^{\prime}}\cdot
\phi_{\Theta}^* c_1(\widehat{E(-\Theta)})>0,$$
where $r_F : = {\rm  rk} (F)$ and $r^{\prime}_E :=
{\rm rk}(\widehat{E(-\Theta)})$.
(Recall that $\phi_{\Theta}$ is the isogeny induced by $\Theta$.)
\end{theorem}

\begin{remark}\label{chern}
Although most conveniently written in terms of the Fourier-Mukai
transform, the statement of the theorem is indeed a numerical
condition intrinsic to $E$ (and $F$), since by \cite{mukai2}
Corollary 1.18 one has:
$$c_1(\widehat{E(-\Theta)})=-PD_{2g-2}({\rm ch}_{g-1}(E(-\Theta))),$$
where $PD$ denotes the Poincar\'e duality map
$$PD_{2g-2}:H^{2g-2}(J(X),\ZZ)\rightarrow H^{2}(J(X),\ZZ),$$
and ${\rm ch}_{g-1}$ the $(g-1)$-st component of the Chern
character. Note also that
$${\rm rk}(\widehat{E(-\Theta)})= h^0(E(-\Theta)) = \frac{1}{r^{g-1}}\cdot
\frac{c_1(E(-\Theta))^g}{g!}$$
by Lemma \ref{strong} and \cite{mukai1} Corollary 2.8.
\end{remark}

We can assume $E$ and $F$ to be simple by the considerations in \S2,
and we will do so in what follows. We begin with a few technical
results. In the first place, it is useful to consider the \emph{skew
Pontrjagin product}, a slight variation of the usual Pontrjagin
product (see \cite{pareschi} \S1). Namely, given two sheaves $\E$
and $\G$ on $X$, one defines
$$\E\hat * \G:={d}_*(p_1^*(\E)\otimes p_2^*(\G)),$$
where $p_1$ and $p_2$ are the projections from $X\times X$ to the
two factors and $d:X\times X\rightarrow X$ is the difference map.

\begin{lemma}\label{mult2}
For all $i\geq 0$ we have:
$$h^i((E\hat{*} F)\otimes \OO_X(-\Theta)) = h^i((E\hat{*}\OO_X(-\Theta))\otimes F).$$
\end{lemma}
\begin{proof}
This follows from Lemma 3.2 in \cite{pareschi} if we prove the
following vanishings:
\begin{enumerate}
\item $h^i(t_x^*E\otimes F) =0,~\forall i>0,~\forall x\in X.$
\item $h^i(t_x^*E\otimes \OO_X(-\Theta))=0,~\forall i>0,~\forall x\in X.$
\end{enumerate}
We treat them separately:
\begin{enumerate}
\item By Lemma \ref{sh-properties}(1) we know that there exist isogenies $\pi_E:Y_E\rightarrow
X$ and $\pi_F: Y_F\rightarrow X$, and line bundles $M$ on $Y_E$ and
$N$ on $Y_F$, such that $\pi_E^*E\cong \underset{r_E}{\oplus} M$ and
$\pi_F^*F\cong \underset{r_F}{\oplus} N$. Now on the fiber product
$Y_E\times_X Y_F$, the pull-back of $t_x^*E\otimes F$ is a direct
sum of line bundles numerically equivalent to $p_1^*M\otimes p_2^*N$
. This line bundle is ample and has sections, and so no higher
cohomology by the Index Theorem. Consequently the same must be true
for $t_x^*E\otimes F$.
\item Since $E$ is semihomogeneous, we have $t_x^*E\cong E\otimes \alpha$ for some
$\alpha\in \Pic0$, and so:
$$h^i(t_x^*E\otimes \OO_X(-\Theta))=h^i(E\otimes\OO_X(-\Theta)\otimes \alpha)=0,$$
since $E(-\Theta)$ satisfies $IT_0$.
\end{enumerate}
\end{proof}

Let us assume from now on for simplicity that the polarization
$\Theta$ is symmetric. This makes the proofs less technical, but the
general case is completely similar since everything depends (via
suitable isogenies) only on numerical equivalence classes.

\begin{proposition}\label{cohomological-multiplication}
Under the hypotheses above, the multiplication maps
$$H^0(E)\otimes H^0(t_x^*F)\longrightarrow H^0(E\otimes t_x^*F)$$
are surjective for all $x\in X$ if we have the following vanishing:
$$h^i(\phi_{\Theta}^*((-1_X)^*E\otimes \OO_X(-\Theta))^{\widehat{}}\otimes F(-\Theta))=0,~\forall i>0.$$
\end{proposition}
\begin{proof}
By \cite{pareschi} Theorem 3.1, all the multiplication maps in the
statement are surjective if the skew-Pontrjagin product $E\hat{*}F$
is globally generated, so in particular if $(E\hat{*} F)$ is
$0$-$\Theta$-regular. On the other hand, by Lemma \ref{mult2}, we
can check this $0$-regularity by checking the vanishing of
$h^i((E\hat{*}\OO_X(-\Theta))\otimes F)$. To this end, we use
Mukai's general Lemma 3.10 in \cite{mukai1} to see that
$$E\hat{*}\OO_X(-\Theta)\cong \phi_{\Theta}^*((-1_X)^*E\otimes \OO_X(-\Theta))^{\widehat{}}
\otimes \OO(-\Theta).$$ This implies the statement.
\end{proof}

We are now in a position to prove Theorem
\ref{mixed-multiplication}: we only need to understand the numerical
assumptions under which the cohomological requirement in Proposition
\ref{cohomological-multiplication} is satisfied.

\medskip
\begin{proof}(\emph{of Theorem \ref{mixed-multiplication}}.)
We first apply Lemma \ref{sh-properties}(1) to
$G:=\phi_{\Theta}^*\widehat{(-1_X)^*E(-\Theta)}$ and
$H:=F(-\Theta)$: there exist isogenies $\pi_G:Y_G\rightarrow X$ and
$\pi_H:Y_H\rightarrow X$, and line bundles $M$ on $Y_G$ and $N$ on
$Y_H$, such that $\pi_G^*G\cong \underset{r_G}{\oplus}M$ and
$\pi_H^*H\cong \underset{r_H}{\oplus}N$. Consider the fiber product
$Z:=Y_G\times_X Y_H$, with projections $p_G$ and $p_H$. Denote by
$p:Z\rightarrow X$ the natural composition. By pulling everything
back to $Z$, we see that
$$p^*(G\otimes H)\cong \underset{r_G\cdot r_F}{\bigoplus}(p_1^*M\otimes p_2^*N).$$
This implies that our desired vanishing $H^i(G\otimes H)=0$ (cf.
Proposition \ref{cohomological-multiplication}) holds as long as
$$H^i(p_G^*M\otimes p_H^*N)=0, ~\forall i>0.$$

Now $c_1(p_G^*M)=p_G^*c_1(M)=\frac{1}{r_G}p^*c_1(G)$ and similarly
$c_1(p_H^*N)=p_H^*c_1(N)=\frac{1}{r_H}p^*c_1(G)$. Finally we get
$$c_1(p_G^*M\otimes p_H^*N)=p^*(\frac{1}{r_G}\cdot c_1(G)+\frac{1}{r_H}\cdot c_1(H)).$$
Thus all we need to have is that the class
$$\frac{1}{r_G}\cdot c_1(G)+\frac{1}{r_H}\cdot c_1(H)$$
be ample, and this is clearly equivalent to the statement of the
theorem.
\end{proof}

\noindent
{\bf $(-1)$-$\Theta$-regular vector bundles.}
It can be easily seen that Theorem \ref{mixed-multiplication} implies that a
$(-1)$-$\Theta$-regular semihomogeneous bundle is normally generated.
Under this regularity hypothesis we have however
a much more general statement, which works for every vector
bundle on a polarized abelian variety.

\begin{theorem}
For $(-1)$-$\Theta$-regular vector bundles $E$ and $F$ on $X$,
the multiplication map
$$H^0(E)\otimes H^0(F)\rightarrow H^0(E\otimes F)$$
is surjective.
\end{theorem}
\begin{proof}
We use an argument exploited in \cite{pp1}, inspired by techniques
introduced by Kempf. Let us consider the diagram
$$\xymatrix{
\bigoplus_{\xi\in U} H^0(E(-2\Theta)\otimes P_\xi)\otimes
H^0(2\Theta\otimes P_\xi^\vee)\otimes H^0(F) \ar[r] \ar[d] &
H^0(E)\otimes
H^0(F) \ar[d]  \\
\bigoplus_{\xi\in U} H^0(E(-2\Theta)\otimes P_\xi)\otimes
H^0(F(2\Theta)\otimes P_\xi^\vee)  \ar[r] & H^0(E\otimes F) }$$
Under the given hypotheses, the bottom horizontal arrow is onto by
the general Theorem \ref{surjectivity}. On the other hand, the
abelian Castelnuovo-Mumford Lemma Theorem \ref{acm} insures that
each one of the components of the vertical map on the left is
surjective. Thus the composition is surjective, which gives the
surjectivity of the vertical map on the right.
\end{proof}

\begin{corollary}
Every  $(-1)$-$\Theta$-regular vector bundle is normally generated.
\end{corollary}

\noindent {\bf Examples.} There are two basic classes of examples of
$(-1)$-$\Theta$-regular bundles, and both turn out to be
semihomogeneous. They correspond to the properties of linear series
on abelian varieties and on moduli spaces of vector bundles on
curves, respectively.

\begin{example}({\bf {Projective normality of line bundles.}})
For every ample divisor $\Theta$ on $X$, the line bundle $L
=\OO_X(m\Theta)$ is $(-1)$-$\Theta$-regular for $m\geq 3$. Thus we
recover the classical fact that $\OO_X(m\Theta)$ is projectively
normal for $m\geq 3$.
\end{example}

\smallskip

\begin{example}({\bf {Verlinde bundles.}})\label{Verlinde}
Let $U_C(r,0)$ be the moduli space of rank $r$ and degree $0$
semistable vector bundles on a a smooth projective curve $C$ of
genus $g\geq 2$. This comes with a natural determinant map ${\rm
det}: U_C(r,0)\rightarrow J(C)$, where $J(C)$ is the Jacobian of
$C$. To a generalized theta divisor $\Theta_N$ on $U_C(r,0)$
(depending on the choice of a line bundle $N\in {\rm Pic}^{g-1}(C)$)
one associates for any $k\geq 1$ the $(r,k)$-Verlinde bundle on
$J(C)$, defined by $E_{r, k}:= {\rm det}_* \OO(k\Theta_N)$ (cf.
\cite{popa}). It is shown in \emph{loc. cit.} that the numerical
properties of $E_{r,k}$ are essential in understanding the linear
series $|k\Theta_N|$ on $U_C(r,0)$. It is noted there that $E_{r,k}$
are polystable and semihomogeneous.

\noindent A basic property of these vector bundles is the fact that
$r_J^*E_{r,k}\cong \oplus \OO_J(kr\Theta_N)$,
where $r_J$ denotes multiplication by $r$ on $J(C)$ (cf. \cite{popa}
Lemma 2.3). Noting that the pull-back $r_J^*\OO_J(\Theta_N)$ is
numerically equivalent to $\OO(r^2\Theta_N)$, we obtain that
$E_{r,k}$ is $0$-$\Theta$-regular iff $k\geq r+1$, and
$(-1)$-$\Theta$-regular iff $k\geq 2r+1$. This implies by the
statements above that $E_{r,k}$ is globally generated for $k\geq
r+1$ and normally generated for $k\geq 2r+1$. These are precisely
the results \cite{popa} Proposition 5.2 and Theorem 5.9(a), the
second obtained there by ad-hoc (though related) methods.
\end{example}

\providecommand{\bysame}{\leavevmode\hbox
to3em{\hrulefill}\thinspace}

\end{document}